\documentclass[12pt]{article}

\usepackage{amsthm,amsmath,amssymb}

\usepackage{graphicx}
 \usepackage{epsfig}
  \usepackage{latexsym, epsfig, psfrag,eepic,colordvi,bm}
  \usepackage{graphics,color}
  \usepackage{graphics}
  \usepackage{graphicx}
  \usepackage{epsfig}
  \usepackage[all]{xy}
\usepackage{cancel}
  \usepackage{amssymb}
  \usepackage{amsthm,amsmath}

  \usepackage{latexsym, epsfig, psfrag,eepic,colordvi,bm}
  \usepackage{graphics,color}
  \usepackage{amsmath,amsfonts,amssymb,amscd}
  \usepackage[all]{xy}
  \usepackage{epstopdf}

  \usepackage{amsthm,amsmath,amssymb}

  \usepackage{graphicx}
  \usepackage{amsthm,amsmath,amssymb}

  \usepackage{graphicx,epsfig}
\usepackage[colorlinks=true,citecolor=black,linkcolor=black,urlcolor=blue]{hyperref}

\usepackage[all]{xy}
\usepackage{graphicx}
\usepackage{mathptmx}

\usepackage{latexsym, epsfig,cite, psfrag,eepic,colordvi,bm}
\usepackage{graphics,color,graphicx}
\usepackage[justification=centering]{caption}
\usepackage{amsmath,amsfonts,amssymb,amscd}
\usepackage[all]{xy}

\theoremstyle{plain}
\newtheorem{theorem}{Theorem}[section]

\newtheorem{lemma}{Lemma}[section]
\newtheorem{corollary}{Corollary}[section]
\newtheorem{proposition}{Proposition}[section]

\theoremstyle{definition}

\theoremstyle{remark}
\newtheorem{remark}[theorem]{Remark}

\newcommand{\genstirlingI}[3]{%
	\genfrac{[}{]}{0pt}{#1}{#2}{#3}%
}
\newcommand{\genstirlingII}[3]{%
	\genfrac{\{}{\}}{0pt}{#1}{#2}{#3}%
}
\newcommand{\stirlingI}[2]{\genstirlingI{}{#1}{#2}}

\newcommand{\stirlingII}[2]{\genstirlingII{}{#1}{#2}}

\date{}

\title{Towards studying the structure of \\triple Hurwitz numbers}

\author{Ricky Xiao-Feng Chen~\footnote{ORCID: 0000-0003-1061-3049}\\
	\small School of Mathematics, Hefei University of Technology \\[-0.8ex]
	\small 485 Danxia Road, Hefei, Anhui Province 230601, P.~R.~China\\[-0.8ex]
	\small\tt xiaofengchen@hfut.edu.cn
}

\begin{document}

\maketitle

\begin{abstract}
	
	Going beyond the studies of single and double Hurwitz numbers, we report some progress towards studying Hurwitz numbers which correspond to ramified coverings of the Riemann sphere involving three nonsimple branch points.
	We first prove a recursion which implies a fundamental identity of Frobenius enumerating factorizations of a permutation
	in group
	algebra theory.
	We next apply the recursion to study Hurwitz numbers involving three nonsimple branch points (besides simple ones),
	two of them having deterministic ramification profiles while the remaining one having a prescribed number of preimages.
	The recursion allows us to obtain recurrences as well as explicit formulas for these numbers which
	also generalize a number of existing results on single and double Hurwitz numbers.

	The case where one of the nonsimple branch points with deterministic profile has a unique preimage (one-part
	quasi-triple Hurwitz numbers {\color{blue} or $(1,m)$-part triple Hurwitz numbers}) is particularly studied in detail.
	We prove an attractive dimension-reduction formula from which any one-part quasi-triple Hurwitz number
	can be reduced to quasi-triple Hurwitz numbers where both with deterministic profiles are fully ramified.
	We also obtain the polynomiality of one-part quasi-triple Hurwitz numbers analogous to that implied by
	the remarkable ELSV formula for single Hurwitz numbers, and discuss the potential connection to intersection theory.


\end{abstract}

\section{Introduction}\label{sec1}

\subsection{Single Hurwitz numbers}

Hurwitz numbers count the non-equivalent ramified coverings of a Riemann surface $M$ of genus $h$
by a Riemann surface $M'$ of genus $g$, which was introduced by Hurwitz in 1891~\cite{hurwitz}. Here we only focus on
the case $h=0$, i.e., the Riemann sphere $\mathbb{S}^2$ (or $\mathbb{CP}^1$).
A ramified covering of degree $d$ from $M$ to $\mathbb{S}^2$ is a nonconstant meromorphic function $f: M\rightarrow \mathbb{S}^2$
such that all but a finite number of points in the sphere have $d$ distinct preimages.
These finite number of points are called branch points of the covering.
The preimages of a branch point are called poles.
Suppose a branch point $p$ has $n$ preimages $z_1, z_2, \ldots, z_n$ and the point $z_i$ has a multiplicity of $\alpha_i$.
Then, $$\alpha_1+\alpha_2+\cdots +\alpha_n=d,$$
 and $z_i$ is called a pole of order $\alpha_i$ while $p$
is called a branch point of ramification profile (or type) 
$$
\alpha=(\alpha_1, \alpha_2, \ldots, \alpha_n).
$$
A branch point with $d-1$ preimages, i.e., $n=d-1$, is called a simple branch point.
Two ramified coverings $f$ and $f'$ are viewed equivalent if there exists a homomophism $h: M\rightarrow M$
such that $f=f'\circ h$.

A formula for the numbers (referred to as single Hurwitz numbers later) corresponding to a single nonsimple branch point besides simple ones and $g=0$ was given by
Hurwitz himself without a complete proof, and rediscovered by Goulden and Jackson~\cite{GJ97}.
The approach of Goulden and Jackson is based on permutation factorization, that is, the desired number equals
the number of ordered factorizations of a permutation of $d$ elements and cycle-type $\alpha$ into
$t$ transpositions for some number $t$ where the transpositions act transitively on the $d$ elements. This approach is presumably close to
the original Hurwitz way of reasoning as can be seen in the reconstruction of Hurwitz's
proof in~\cite{recon96}.
Subsequently, explicit formulas of these numbers $\mathsf{H}^g_{\alpha}$ for genus $g=1$ and $g=2$ were given in
Goulden and Jackson~\cite{GJ99,GJ99a}, Goulden, Jackson and Vainshtein~\cite{GJV00}.
A number of interesting properties like polynomiality~\cite{GJV00} have been studied too.
In addition, Li, Zhao and Zheng~\cite{lzz00}
obtained a recurrence for single Hurwitz numbers via the relative Gromov-Witten theory.

\subsection{Connections to Hodge integrals}
The single Hurwitz numbers $\mathsf{H}^g_{\alpha}$ have also been found very closely connected with the geometry of the moduli spaces of curves,
with diverse physical theories, and with singularity theory. In particular, there are remarkable connections with the Hodge integrals firstly exhibited by the celebrated ELSV formula~\cite{ELSV1,ELSV}:
\begin{align}
\mathsf{H}^g_{\alpha}= t! \prod_{i=1}^n \frac{\alpha_i^{\alpha_i}}{\alpha_i !} \int_{\overline{{M}}_{g,n}} \frac{1-\lambda_1+\lambda_2-
\cdots \pm \lambda_g}{(1-\alpha_1 \psi_1)\cdots (1-\alpha_n \psi_n)},
\end{align}
where $\overline{{M}}_{g,n}$ is the moduli space of stable curves of genus $g$ with $n$ marked points,
$\lambda_i$ is a certain codimension $i$ class, and $\psi_k$ is a certain
codimension $1$ class (Chern class).
We refer to~\cite{ELSV1,ELSV} for precise definitions of these objects.
Properties, such as, the quantity
$$\mathsf{H}^g_{\alpha} \Big/ t! \prod_{i=1}^n \frac{\alpha_i^{\alpha_i}}{\alpha_i !}$$
is a polynomial in $\alpha_i$'s succinctly follow from the ELSV formula, settling
a conjecture in~\cite{GJV00}.
With this, the intersection numbers
$$
\int_{\overline{{M}}_{g,n}} \psi_1^{b_1} \psi_2^{b_2} \cdots \psi_n^{b_n} \lambda_k
$$
are just the coefficients of the polynomial, and
 many problems in intersection theory of curves, e.g., the $\lambda_g$ conjecture~\cite{FP},
may be solved by studying Hurwitz numbers.

Along these lines, other Hurwitz numbers such as double Hurwitz numbers~\cite{GJV05,DL22,wedge15}, double Hurwitz numbers with completed cycles~\cite{3-cycle,complete11}, monotone
Hurwitz numbers~\cite{monotone} and orbifold Hurwitz numbers~\cite{orbiford} etc.~have been extensively investigated, and finding other connections such as searching their Hodge integral counterparts has attracted
lots of attention as well. These existing Hurwitz numbers have been also studied in the framework of KP hierarchy, see e.g.~\cite{GJ08,bdks,kaz-lan15} (in particular,
\cite{bdks} may cover all existing Hurwitz numbers). 

\subsection{Our objective: studying triple Hurwitz numbers}

 Single and double Hurwitz numbers are defined according to the (maximum) number of involved nonsimple branch points
 in the considered ramified coverings of the Riemann sphere: single for one nonsimple branch point and
 double for two nonsimple branch points.
 It is worth noting that the existing results for terms without the word single or double, e.g., monotone
 and orbifold Hurwitz numbers, can be also viewed as just for certain classes of single and double Hurwitz numbers.
 
Our objective in the present paper is to study genus $g$ ramified coverings of $\mathbb{S}^2$ with up to three nonsimple branch points. We call the numbers of these coverings triple Hurwitz numbers
 in analogy to single and double Hurwitz numbers.
 In particular, we compute the number corresponding to the degree $d$ coverings (possibly disconnected) of genus $g$,
 where besides simple branch points, the ramification profile of the first nonsimple branch point (say $0$) is $\alpha$, the ramification profile of the second nonsimple branch point (say $\infty$)
 is $\beta=(\beta_1,\beta_2,\ldots, \beta_n)$,
 and the ramification profile of the third nonsimple branch point is any one of the form $(\gamma_1,\gamma_2,\ldots, \gamma_m)$ for a fixed $m$.
 We call these numbers quasi-triple Hurwitz numbers to distinguish them from ``standard" triple Hurwitz
 numbers involving three nonsimple branch points all with deterministic ramification profiles.
 We denote the quasi-triple Hurwitz number by $\mathsf{H}_{d,m}^{g}(\alpha, \beta)$.
 
 The first general explicit formula for single Hurwitz numbers appeared over 100 years ago (i.e., Hurwitz's formula for genus zero), while a
 systematic study of double Hurwitz numbers also has a history of more than 20 years.
 However, to the best of our knowledge, no explicit formulas and properties have been obtained
 for Hurwitz numbers involving three or more nonsimple branch points in addition to 
 simple branch points.
 This is of course not because these are not interesting but intrinsically hard to solve, and
 here we attempt to make some progress on studying the structure of these triple Hurwitz numbers.

Following Goulden and Jackson and others, our approach for tackling triple Hurwitz numbers is also based on enumerating some permutation products.
The enumeration of permutation products or factorizations subject to certain conditions in the symmetric groups
has a long rich history. Its importance lies in that many problems arising in other areas, such as geometry,
mathematical physics and combinatorics,
can be translated into algebraic formulation of permutation products or factorizations.
For instance, besides Hurwitz numbers, graph embeddings in topological surfaces (or maps)
can be formulated into permutation products subject to the Euler characteristic formula~\cite{edmonds,jac,walsh1,walsh2}, biological structures
like RNA structures and genome sequence rearrangement problems~\cite{aprw,chr1} also have permutation product representations
and solutions.

The core of our approach is from symmetric group characters and symmetric functions which
have already been explored with great success,
in particular, in a series of work of Jackson, Goulden and their coauthors~\cite{jac,jac88,GJ97,GJ99,GJ99a,GJV00}.
However, we shall see 
our approach surprisingly requires a little knowledge of the character theory
and symmetric function theory.

The main results of the present paper are: 
\begin{itemize}
\item A general symmetric function description of $\mathsf{H}_{d,m}^{g}(\alpha, \beta)$;
\item Explicit formulas for
one-part quasi-triple Hurwitz numbers (or $(1,m)$-part triple Hurwitz numbers) $\mathsf{H}_{d,m}^{g}((d), \beta)$, i.e., these with $\alpha$ having only one part;
\item Based on the explicit formulas, we derive: (i)
an efficient recurrence for computing one-part quasi-triple Hurwitz numbers of all genera and profiles,
and (ii) the polynomiality of one-part quasi-triple Hurwitz numbers and possible connections to intersection theory of curves.
\end{itemize}
Besides, that our computations of triple Hurwitz numbers become feasible is due to some new results on character theory of symmetric groups which
are of independent interest.
Note that $\mathsf{H}_{d,d}^{g}(\alpha, \beta)$, i.e., $m=d$, are double Hurwitz numbers.
As such, our results recover and generalize a number of existing results on single and double Hurwitz numbers as well.

\subsection{Plan of the paper}
The organization of the paper is as follows.
In Section~\ref{sec2}, we introduce relevant symmetric
group character theory and prove a novel recursion
about permutation products.
The novel recursion generalizes a fundamental identity of Frobenius in group
character theory in the context of symmetric groups. A direct, compact enumeration formula (Theorem~\ref{thm:recur2}) is obtained as well.
The proof of the latter
is based upon studying
the M\"{o}bius functions of partially ordered sets consisting of 
set partitions, which is separately presented in Section $4$.
These are the main ingredients of our approach.

In Section~\ref{sec3}, we apply the novel recursion and compact formula to the enumeration of quasi-triple
Hurwitz numbers.
We first provide a general symmetric function description (Theorem~\ref{thm:shur-express}) making
use of the power sum symmetric functions and Schur symmetric functions.

Next, we focus on one-part quasi-triple Hurwitz numbers (i.e., one of the three nonsimple
branch points is fully ramified).
We obtain some recurrences and explicit formulas for these numbers.
The involved quantities are expressed in terms of the hyperbolic $\sinh$ function (Theorem~\ref{thm:triple-sinh})
and the generalized Bernoulli polynomials (Theorem~\ref{thm:triple}).
	We also prove a dimension-reduction formula from which any one-part quasi-triple Hurwitz number
can be reduced to quasi-triple Hurwitz numbers where two branch points are fully ramified.
This contribution advances the state of the art that Hurwitz numbers with at most two nonsimple
branch points can be explicitly computed.

Based on an analysis of the obtained formulas, we subsequently prove a polynomiality property (Theorem~\ref{thm:coeffi}) of the one-part
quasi-triple Hurwitz numbers in analogy with the single and double Hurwitz numbers.
As discussed throughly, the polynomiality suggests a potential connection to some ELSV-type formulas
concerning intersection theory over certain space of curves.
Accordingly, we discuss these speculated intersection numbers.
In particular, as for one-part double Hurwitz numbers, we are able to present
a clean formula for all intersection numbers there.
These results may facilitate the search of these ELSV-type formulas in the future.

\section{A novel recursion}\label{sec2}

Let $[d]:=\{1,2,\ldots,d\}$, and denote by $\mathfrak{S}_d$ the symmetric group on $[d]$,
i.e., the set of permutations on $[d]$ with the functional composition as the multiplication.
Let $\pi \in \mathfrak{S}_d$ and $x\in [d]$.
There obviously exists a number $p\leq d$ such that $x, \pi(x), \pi^2(x), \ldots, \pi^{p-1}(x)
$ are all distinct while $\pi^{p}(x)=x$. This gives a cycle of $\pi$ and $p$ is called its length.
A permutation can be written as a product of its disjoint cycles generated in this manner.
The length distribution of the cycles of $\pi$ is called the cycle-type of $\pi$ and can be encoded as an (integer) partition
of $d$. A partition $\lambda$ of $d$, denoted by $\lambda \vdash d$,
is usually represented by a nonincreasing positive integer sequence $\lambda=(\lambda_1, \lambda_2,\ldots, \lambda_k)$
such that $\lambda_1+\lambda_2+\cdots+\lambda_k=d$. The number $k$ is called the length of $\lambda$, denoted by $\ell(\lambda)$.
Another representation of a partition is in the form $\lambda=[1^{m_1}, 2^{m_2},\ldots, d^{m_d}]$
indicating that there are $m_i$ of $i$'s in the partition.
We usually discard the entry $i^{m_i}$ for $m_i=0$ and write the entry $i^1$ simply as $i$.
Clearly, $\sum_{i=1}^d i m_i=d$ and $\ell(\lambda)=\sum_{i=1}^d m_i$.
We also define $Aut(\lambda)=\prod_{i=1}^d m_i!$, that is, the number of permutations of the entries of $\lambda$ that
fix $\lambda$. 
We will use the two representations of partitions interchangably, whichever is more convenient.

It is well known that a conjugacy class of $\mathfrak{S}_d$ consists of
permutations of the same cycle-type.
We denote the conjugacy class indexed by $\lambda$ by $\mathcal{C}_{\lambda}$.
If $\pi \in \mathcal{C}_{\lambda}$, then the number of cycles contained in $\pi$ is $\ell(\pi)=\ell(\lambda)$.
Moreover, the number of elements contained in $\mathcal{C}_{\lambda}$ is well known to be $|\mathcal{C}_{\lambda}|=\frac{d!}{z_{\lambda}}$
where
$$
z_{\lambda}= \prod_{i=1}^d i^{m_i} m_i!.
$$
And the number of permutations having exactly $k$ cycles is given by the signless Stirling
number of the first kind $\stirlingI{d}{k}$ which can be characterized by:
$$
\sum_{k=1}^d (-1)^{d-k}\stirlingI{d}{k} x^k = x(x-1)\cdots (x-d+1).
$$

We write the character associated to the irreducible representation indexed by $\lambda$
as $\chi^{\lambda}$. Suppose $C$ is a conjugacy class of $\mathfrak{S}_d$ indexed by $\alpha$ and $\pi \in C$.
We view $\chi^{\lambda}(\pi), \, \chi^{\lambda}(\alpha), \, \chi^{\lambda}(C)$ as the same, and we trust the context to prevent confusion.

The next is basically a fundamental result due to Frobenius (some attribute it to Burnside)
in the group algebra theory, in particular, the center of the group algebra,
and we just reformulate it in an unconventional way.

\begin{theorem}[Frobenius identity]\label{thm:fro}
	Let $C_i$, for $1\leq i \leq k$, be a conjugacy class of $\mathfrak{S}_d$ and $N_{C_1,C_2,\ldots, C_k}(\pi)$ denote the number of tuples $(\sigma_1, \sigma_2,\ldots, \sigma_k)$
	such that $\sigma_i \in C_i$ and $\sigma_1 \sigma_2 \cdots \sigma_k=\pi$.
	Suppose the function $N_{C_1,C_2,\ldots, C_k}=\sum_{\lambda \vdash d} \mathfrak{n}_{\lambda} \chi^{\lambda}$.
	Then,
	\begin{align}
\mathfrak{n}_{\lambda}=\frac{\prod_{i=1}^{k}|C_i|}{d!}  \{ \dim(\lambda) \}^{-k+1} \prod_{i=1}^k \chi^{\lambda}(C_i),
	\end{align}
where $\dim(\lambda)$ stands for the dimension of the irreducible representation indexed by $\lambda$.
\end{theorem}

The forthcoming results are well known and also attributed to Frobenius:
\begin{align}
	\eta(\lambda)=\frac{|\mathcal{C}_{[1^{d-2}, \, 2]}| \, \chi^{\lambda}([1^{d-2}, \, 2])}{\dim(\lambda)}= \sum_i {\lambda_i \choose 2} -\sum_i {\widetilde{\lambda}_i \choose 2},
\end{align}
where $\widetilde{\lambda}$ stands for the conjugate of $\lambda$.

Now we are ready to present a novel recursion regarding permutation products.
Let 
\begin{align*}
	\mathfrak{c}_{\lambda,m} &= \sum_{j=0}^m (-1)^j {m \choose j} \mathfrak{m}_{\lambda, m-j}, \quad 	\mathfrak{m}_{\lambda,m}=\prod_{u\in \lambda}\frac{m+c(u)}{h(u)},
\end{align*}
where for $u=(i,j)\in \lambda$, i.e., the $(i,j)$ cell in the Young diagram of $\lambda$, $c(u)=j-i$ and $h(u)$ is the hook length of 
the cell $u$.

\begin{theorem}\label{thm:recur2}
	Let $\xi_{d,m}(C_1,\ldots, C_t)$ be the number of tuples $( \sigma_1,\sigma_2,\ldots ,\sigma_t)$ 
	such that the permutation $\pi=  \sigma_1\sigma_2\cdots \sigma_t$ has $m$ cycles, where $\sigma_i $
	belongs to a conjugacy class ${C}_i$ of $\mathfrak{S}_d$. Then, we have
	\begin{align}
		\xi_{d,m} 
		&=\sum_{k\geq 0}    (-1)^k \stirlingI{m+k}{m} W_{d,m+k}, \label{eq:compact-direct}
	\end{align}
	where
\begin{align}\label{eq:WW}
	W_{d,m}(C_1,\ldots, C_t)= \frac{\prod_{i=1}^{t}|C_i|}{m!} \sum_{\lambda \vdash d} \mathfrak{c}_{\lambda,m}	 \big\{\dim(\lambda)\big\}^{-t+1} \prod_{i=1}^t \chi^{\lambda}(C_i).
\end{align}

\end{theorem}
We remark that Theorem~\ref{thm:recur2} can be further generalized where the permutations $\pi$ are not necessarily from products of permutations in prescribed conjugacy classes and we will report the progress on this in the future.

\section{Quasi-triple Hurwitz numbers}\label{sec3}

Recall $\mathsf{H}_{d,m}^{g}(\alpha, \beta)$ is the triple Hurwitz number corresponding to the genus $g$
ramified coverings of $\mathbb{S}^2$ with two nonsimple branch points of profiles respectively
$\alpha$ and $\beta$ and with one nonsimple branch point having $m$ preimages.
According to the Riemann-Hurwitz theorem, the number $t$ of simple branch points in these coverings satisfies:
\begin{align} \label{eq:R-H}
t=2g-2+m-d+ \ell(\beta)+\ell(\alpha).
\end{align}
It is noted that both the cases $m=d$ and $m=d-1$ reduce to certain double Hurwitz numbers~\cite{GJV05}.
In the following,
we will first prove a symmetric function description of the general quasi-triple 
Hurwitz numbers and then focus on one-part quasi-triple Hurwitz numbers, i.e., those where
$\alpha=(d)$ and thus $\ell(\alpha)=1$.

\subsection{A general symmetric function description}

It is sometimes more convenient to employ the following notation relating to quasi-triple Hurwitz numbers:
$$
\mathfrak{H}_{d,m}^{t}(\alpha, \beta) = \xi_{d,m}(\alpha,\mathcal{C}_{[2, 1^{d-2}]}, \ldots, \mathcal{C}_{[2, 1^{d-2}]}, \beta ),
$$
where there are $t \geq 0$ of $\mathcal{C}_{[2, 1^{d-2}]}$.
In this notation, the number $t$ is not restricted to the Riemann-Hurwitz theorem.
Define
$$
\mathfrak{H}_{d}=\sum_{m>0} \sum_{t\geq 0} \sum_{\alpha \vdash d, \, \beta \vdash d} \mathfrak{H}_{d,m}^{t}(\alpha, \beta) p_{\alpha}({\bf x}) p_{\beta}({\bf y}) \frac{z^t}{t!} w^m,
$$
where $p_{\lambda}({\bf x})$ stands for the power sum symmetric function indexed by $\lambda$ in indeterminates ${\bf x}=(x_1, x_2,\ldots)$.
For the definitions of power sum symmetric functions $p_{\lambda}$ and Schur symmetric functions $s_{\lambda}$ used later, we follow Stanley~\cite{ec-2}.

With this, we are able to prove the following symmetric function characterization.

\begin{theorem}\label{thm:shur-express}
	 For $d >0$, we have
	\begin{align}
	\frac{\mathfrak{H}_{d}}{(d!)^2}
	=	&  \sum_{\lambda \vdash d} \,
	\frac{\mathfrak{m}_{\lambda, w} }{  \dim(\lambda)}  e^{\eta(\lambda)z} { s_{\lambda}({\bf x})s_{\lambda}({\bf y})}.
\end{align}
\end{theorem}

\proof According to Theorem~\ref{thm:recur2}, for fixed $t\geq 0, \alpha, \beta$, we compute the corresponding $W$-numbers:
\begin{align*}
	W_{d,m+k} &= \sum_{\lambda \vdash d} \mathfrak{c}_{\lambda,m+k}	\frac{|\mathcal{C}_{\alpha}| |\mathcal{C}_{\beta}| |\mathcal{C}_{[1^{d-2}, 2]}|^t}{(m+k)!}  \big\{\dim(\lambda) \big\}^{-t-1}  \{ \chi^{\lambda}([1^{d-2}, 2]) \}^t \chi^{\lambda}(\alpha) \chi^{\lambda}(\beta)\\
	&=\sum_{\lambda \vdash d} \mathfrak{c}_{\lambda,m+k}	\frac{|\mathcal{C}_{\alpha}| |\mathcal{C}_{\beta}| \{\eta(\lambda) \}^t }{(m+k)! \, \dim(\lambda)}   \chi^{\lambda}(\alpha) \chi^{\lambda}(\beta).
\end{align*} 
Then, we have
\begin{align*}
	&\mathfrak{H}_{d,m}=	\sum_{t\geq 0} \sum_{\alpha \vdash d, \, \beta \vdash d} \mathfrak{H}_{d,m}^{t}(\alpha, \beta) p_{\alpha}({\bf x}) p_{\beta}({\bf y}) \frac{z^t}{t!}\\
	=& \sum_{t\geq 0} \sum_{\alpha \vdash d, \, \beta \vdash d} \sum_{k\geq 0} (-1)^k \stirlingI{m+k}{m} \sum_{\lambda \vdash d} \mathfrak{c}_{\lambda,m+k}	\frac{|\mathcal{C}_{\alpha}| |\mathcal{C}_{\beta}| \{\eta(\lambda) \}^t }{(m+k)! \, \dim(\lambda)}   \chi^{\lambda}(\alpha) \chi^{\lambda}(\beta) p_{\alpha}({\bf x}) p_{\beta}({\bf y}) \frac{z^t}{t!}\\
		=	&   \sum_{\lambda \vdash d} \sum_{k\geq 0}	\frac{ (-1)^k \stirlingI{m+k}{m} d!^2 }{(m+k)! }
	 \frac{\mathfrak{c}_{\lambda,m+k}   }{  \dim(\lambda)}  e^{\eta(\lambda)z} { s_{\lambda}({\bf x})s_{\lambda}({\bf y})}.
\end{align*}
In the computation above, we used the well-known fact that
$$
\frac{1}{d!} \sum_{\pi \in \mathfrak{S}_d} \chi^{\lambda}(\pi) p_{ct(\pi)}({\bf x})=s_{\lambda}({\bf x}),
$$
where $ct(\pi)$ denotes the cycle-type of $\pi$.
Note that $\stirlingI{i}{m}=0$ if $i<m$ and $\mathfrak{c}_{\lambda, i}=0$ if $i>d$.
We proceed to compute
\begin{align*}
	\sum_{m>0} \mathfrak{H}_{d,m} w^m &=    \sum_{\lambda \vdash d} \sum_{i \geq 1} \sum_{m>0}	\frac{ (-1)^{i-m} \stirlingI{i}{m} d!^2 }{(i)! }
	\frac{\mathfrak{c}_{\lambda,i}   }{  \dim(\lambda)}  e^{\eta(\lambda)z} { s_{\lambda}({\bf x})s_{\lambda}({\bf y})} w^m\\
	&= d!^2 \sum_{\lambda \vdash d} \sum_{i \geq 1} {w \choose i}	
	\frac{\mathfrak{c}_{\lambda,i}   }{  \dim(\lambda)}  e^{\eta(\lambda)z} { s_{\lambda}({\bf x})s_{\lambda}({\bf y})} \\
		& =d!^2 \sum_{\lambda \vdash d} 
	\frac{\mathfrak{m}_{\lambda, w} }{  \dim(\lambda)}  e^{\eta(\lambda)z} { s_{\lambda}({\bf x})s_{\lambda}({\bf y})}.
\end{align*}
The last formula follows from the tableau interpretation of $\mathfrak{c}_{\lambda,i} $ and $\mathfrak{m}_{\lambda, i}$,
and the proof is complete.
\qed

Derivation of information about $\mathfrak{H}_{d,m}^{t}(\alpha, \beta)$ from Theorem~\ref{thm:shur-express}
is left for future exploration,
e.g., the connections to KdV or KP hierarchy.

\subsection{One-part quasi-triple Hurwitz numbers}

Similar to the single~\cite{GJ97} and double Hurwitz numbers~\cite{GJV05}, using the permutation model (or called Hurwitz's axioms~\cite{GJV05}),
the one-part quasi-triple Hurwitz number is
\begin{align}\label{eq:hur-def}
\mathsf{H}_{d,m}^{g}(\beta):=\mathsf{H}_{d,m}^{g}((d),\beta)=\frac{Aut((d)) Aut(\beta)}{d!} \xi_{d,m}(\mathcal{C}_{(d)}, C_1,\ldots,C_t,C_{t+1}),
\end{align}
where $C_i=\mathcal{C}_{[1^{d-2}, 2]}$ for $1 \leq i \leq t$, $C_{t+1}=\mathcal{C}_{(\beta_1, \beta_2, \ldots, \beta_n)}$,
and $t=2g-1+n+m-d$.
Here we treat the branch points at $0$ and $\infty$ as labelled whence the weight $\frac{Aut((d)) Aut(\beta)}{d!}$.

\begin{remark}
It is true that when applying the Frobenius identity, the summation over all partitions can be greatly simplified if one of the involved
conjugacy classes is $\mathcal{C}_{(d)}$, i.e., only hook-shapes matter.
As such, one-part single and double Hurwitz numbers can be computed by directly applying the Frobenius identity, e.g., Goulden, Jackson
and Vakil~\cite{GJV05}.
However, one will immediately run into difficulties when there are more than one general conjugacy class involved.
Therefore, it is not that easy to obtain formulas for one-part quasi-triple Hurwitz numbers as one may speculate just
because they are the ``one-part" family.
\end{remark}

To present our result, we first introduce the following numbers and functions.
The hyperbolic sin function $\sinh(x)=(e^x-e^{-x})/2$.
The Bernoulli polynomials of order $N$, ${\bf B}_n^{(N)}(t)$, is defined by the generating function as follows:
$$
B^{(N)}(x,t)=\left(\frac{x}{e^x-1}\right)^N e^{xt}=\sum_{n\geq 0} {\bf B}_n^{(N)}(t) \frac{x^n}{n!}.
$$
The number $B_n^{(N)}:={\bf B}_n^{(N)}(0)$ is called the $n$-th Bernoulli number
of order $N$.
The case $N=1$ gives the classical Bernoulli polynomials $B_n(t)$ and Bernoulli numbers $B_n$.
The following identities hold (see e.g. Srivastava and Todorov~\cite{bernoulli}):
\begin{align}
	{\bf B}_n^{(N)}(x) & =\sum_{k=0}^n {n\choose k} B_k^{(N)} x^{n-k} =\sum_{k=0}^n {n\choose k} B_{n-k}^{(N)} x^{k},\\
	B_n^{(N)}&=\sum_{k=0}^n (-1)^k {N+n \choose n-k}{N+k-1 \choose k}{n+k \choose k}^{-1} S(n+k, k).
	\end{align}
Moreover, the notation $[x^n]f(x)$ means the coefficient of the term $x^n$
in the power series expansion of the function $f(x)$.
First, we have 
	\begin{align}\label{eq:W}
	W_{d,m}(\mathcal{C}_{(d)}, C_1,\ldots, C_t)=\frac{(d-1)!\prod_{i=1}^{t}|C_i|}{m!}\sum_{j=0}^{d-1} (-1)^j \frac{{d-1-j \choose d-m} }{{d-1\choose j}^{t-1}} \prod_{i=1}^t \chi^{[1^j,d-j]}(C_i).
\end{align}

\begin{theorem}[sinh-function formula]\label{thm:triple-sinh}
	Let $g\geq 0$, $m>0$, $d>0$, $\beta=(\beta_1,\ldots, \beta_n)=[1^{a_1},\ldots, d^{a_d}] \vdash d$ and $t=2g+n-1+m-d$.
	Suppose $c_1=a_1-1$ and $c_i=a_i$ for $1<i \leq d$.
	Then,
	the one-part quasi-triple Hurwitz numbers are given by
\begin{align}\label{eq:H-diff}
		\mathsf{H}_{d,m}^{g}(\beta)&={ t! d^{t-1}} [{y^{2g+m-d}}] 
		 \left. \bigg\{ \frac{\partial^d}{ \partial x^d} \, x^{\frac{d-1}{2}} \cdot \frac{ (\log x)^m }{m!}  \cdot \prod_{i=1}^d \Big(\frac{\sinh(\frac{i}{2}\log x+\frac{i}{2}y)}{iy/2} \Big)^{c_i}
	 \bigg\} \right|_{x=1} .
	\end{align}

\end{theorem}

\proof 
In this setting, eq.~\eqref{eq:compact-direct} reduces to
$$
\mathsf{H}_{d,m}^{g}(\beta)=\frac{Aut(\beta)}{d!} \mathfrak{H}_{d,m}^{t}((d), \beta)= \sum_{k=0}^{d-m} (-1)^k \stirlingI{m+k}{m}\widetilde{W_{d,m+k}},
$$
where $\widetilde{W_{d,m+k}}=\frac{Aut(\beta)}{d!} W_{d,m+k}$. According to eq.~\eqref{eq:W},
\begin{align*}
		\widetilde{W_{d,m}}
	&=\frac{(d-1)! {d\choose 2}^t }{m!\prod_{i=1}^n \beta_i}\sum_{j=0}^{d-1} {d-1-j \choose d-m}   {d-1\choose j}^{-t} (-1)^j \left\{\chi^{[1^j, d-j]}([1^{d-2}, 2])\right\}^t \chi^{[1^j, d-j]}(\beta) \\
	&=\frac{(d-1)!  d^t}{m!(d-m)! \prod_{i=1}^n \beta_i}\sum_{j=0}^{d-1} (d-1-j)_{d-m}  (-1)^j \left(\frac{(d-1)}{2}-j \right)^t \chi^{[1^j, d-j]}(\beta)\\
	&=\frac{(d-1)!  d^t}{m!(d-m)! \prod_{i=1}^n \beta_i}\sum_{j=0}^{d-1} \left\{\left. \frac{\mathrm{d}^{d-m}}{\mathrm{d}x^{d-m}} x^{d-1-j} \right|_{x=1} \right\} (-1)^j \left\{ \left[\frac{y^t}{t!}\right] e^{\frac{(d-2j-1)}{2}y} \right\} \chi^{[1^j, d-j]}(\beta)\\
	&=\frac{(d-1)!  d^t}{m!(d-m)! \prod_{i=1}^n \beta_i} \left[\frac{y^t}{t!}\right] \left.\frac{\partial^{d-m}}{\partial x^{d-m}}
	\left\{x^{d-1}e^{\frac{(d-1)y}{2}} \sum_{j=0}^{d-1}  x^{-j}   (-1)^j e^{-jy}  \chi^{[1^j, d-j]}(\beta)
	\right\}\right|_{x=1}.
	\end{align*}
Applying Proposition~\ref{prop:hook-poly} to the last formula, we then have
\begin{align*}
	\widetilde{W_{d,m}} &=\frac{(d-1)!  d^t}{m!(d-m)! \prod_{i=1}^n \beta_i} \left[\frac{y^t}{t!}\right] \left.\frac{\partial^{d-m}}{\partial x^{d-m}}
	\left\{x^{d-1}e^{\frac{(d-1)y}{2}} (1-\frac{e^{-y}}{x})^{-1} \prod_{i=1}^d \Big(1-(\frac{e^{-y}}{x})^i \Big)^{a_i}
	\right\}\right|_{x=1}\\
		&=\frac{(d-1)!  d^t}{m!(d-m)! \prod_{i=1}^n \beta_i} \left[\frac{y^t}{t!}\right] \left.\frac{\partial^{d-m}}{\partial x^{d-m}}
	\left\{x^{d-1} \frac{1}{\prod_i (x^{i/2})^{c_i}} \prod_{i=1}^d \Big(x^{i/2} e^{iy/2}-x^{-i/2}{e^{-iy/2}} \Big)^{c_i}
	\right\}\right|_{x=1}\\
		&=\frac{(d-1)!  d^t}{m!(d-m)! \prod_{i=1}^n \beta_i} \left[\frac{y^t}{t!}\right] \left.\frac{\partial^{d-m}}{\partial x^{d-m}}
	\left\{2^{n-1} x^{\frac{d-1}{2}}  \prod_{i=1}^d \Big(\sinh(\frac{i}{2}\log x+\frac{i}{2}y) \Big)^{c_i}
	\right\}\right|_{x=1}\\
		&=\frac{(d-1)!  d^t}{m!(d-m)! } \left[\frac{y^{t+1-n}}{t!}\right] \left.\frac{\partial^{d-m}}{\partial x^{d-m}}
	\left\{ x^{\frac{d-1}{2}}  \prod_{i=1}^d \Big(\frac{\sinh(\frac{i}{2}\log x+\frac{i}{2}y)}{iy/2} \Big)^{c_i}
	\right\}\right|_{x=1}.
\end{align*} 
The formula for $\widetilde{W_{d,m+k}}$ is analogous.
To summarize, we now obtain
\begin{align}\label{eq:hur-recur}
		\mathsf{H}_{d,m}^{g}(\beta)&={ t! d^{t-1}} [{y^{t+1-n}}]  \sum_{k=0}^{d-m} {d \choose m+k} (-1)^k \stirlingI{k+m}{m} \left. \frac{\partial^{d-m-k}}{ \partial x^{d-m-k}}  x^{\frac{d-1}{2}}  \prod_{i=1}^d \Big(\frac{\sinh(\frac{i}{2}\log x+\frac{i}{2}y)}{iy/2} \Big)^{c_i} \right|_{x=1} .
	\end{align}
Next note that 
\begin{align*}
\frac{(\log x)^m}{m!} = \frac{(\log (1+x-1))^m}{m!} =\sum_{k \geq 0} (-1)^k \stirlingI{m+k}{m} \frac{(x-1)^{m+k}}{(m+k)!}.
\end{align*}
Then, we have
\begin{align*}
\frac{\partial^{m+k}}{ \partial x^{m+k}} \frac{(\log x)^m}{m!} \Bigg|_{x=1} = (-1)^k \stirlingI{m+k}{m},
\end{align*}
and $\frac{\partial^{q}}{ \partial x^{q}} \frac{(\log x)^m}{m!} \Big|_{x=1} =0$ if $q<m$.
Consequently,
we can rewrite eq.~\eqref{eq:hur-recur} as
\begin{align*}
		\mathsf{H}_{d,m}^{g}(\beta)
&= { t! d^{t-1}} [{y^{2g+m-d}}] \left.
		 \bigg\{ \frac{\partial^d}{ \partial x^d}  x^{\frac{d-1}{2}} \frac{(\log x )^m }{m!}  \prod_{i=1}^d \Big(\frac{\sinh(\frac{i}{2}\log x+\frac{i}{2}y)}{iy/2} \Big)^{c_i}
		\bigg\}  \right|_{x=1},
	\end{align*}
as desired. \qed

\begin{theorem}[Bernoulli-polynomial formula]\label{thm:triple}
	Let $g\geq 0$, $m>0$, $d>0$, $\beta=(\beta_1,\ldots, \beta_n)=[1^{a_1},\ldots, d^{a_d}] \vdash d$ and $t=2g+n-1+m-d$.
	Suppose $c_1=a_1-1$ and $c_i=a_i$ for $1<i \leq d$.
	Then,
	the one-part quasi-triple Hurwitz numbers are given by
\begin{align}
	\mathsf{H}_{d,m}^{g}(\beta)=\sum_{k=0}^{d-m} (-1)^k \stirlingI{m+k}{m}\widetilde{W_{d,m+k}},
\end{align}
	where
	\begin{align}
		\widetilde{W_{d,m+k}}
		&= \frac{t! (d-1)!  d^t }{(m+k)! \prod_{i=1}^n \beta_i}  	\sum_{j=0}^{d-m-k} \sum_{1\leq k_1<\cdots <k_p \leq n, \atop p\geq 0}  \frac{(-1)^{m+k-d-p}}{(t+j+1)!}
		\nonumber \\
		&\qquad \qquad \qquad \times {\beta_{k_1}+\cdots +\beta_{k_p}\choose d-m-k-j}
		{\bf B}_{t+j+1}^{(j+1)}\bigg(\frac{1-d}{2}+\sum_{q\notin\{k_1,\ldots, k_p\}} \beta_q \bigg) \label{eq:H-W}.
	\end{align}
	
\end{theorem}

\subsection{Polynomiality}

As a consequence of the ELSV formula~\cite{ELSV}, the Witten symbols, $\langle \cdots \rangle_g$, are connected to the coefficients of the single
Hurwitz numbers:
\begin{align}
	\langle \tau_{b_1}\cdots \tau_{b_n} \lambda_k \rangle_g:= \int_{\overline{\mathcal{M}}_{g,n}} \psi_1^{b_1} \ldots \psi_n^{b_n} \lambda_k=(-1)^k
	[\beta_1^{b_1}\cdots \beta_n^{b_n}] \mathsf{H}^g_{\beta}/t! \prod_{i=1}^n \frac{\beta_i^{\beta_i}}{\beta_i !}.
\end{align}
Specially, the $\lambda_g$-conjecture~\cite{FP} is concerned with an explicit formula for $\langle \tau_{b_1}\cdots \tau_{b_n} \lambda_g \rangle_g$
where $b_1+\cdots+b_n=2g-3+n$:
\begin{align}
	\langle \tau_{b_1}\cdots \tau_{b_n} \lambda_g \rangle_g =\int_{\overline{\mathcal{M}}_{g,n}} \psi_1^{b_1} \ldots \psi_n^{b_n} \lambda_g={2g-3+n \choose b_1,\ldots, b_n} \frac{2^{2g-1}-1}{2^{2g-1}(2g)!} |B_{2g}|.
\end{align}
Thus, the connection between the Hurwitz numbers and Hodge integrals over space of curves can help understand
each other erea better.
In this spirit, researchers are interested in finding Hodge integral (ELSV-type) expressions of various Hurwitz
numbers. For example, for double Hurwitz numbers, there was a conjectured ELSV-type expression in Goulden, Jackson and Vakil~\cite{GJV05}, and
Do and Lewa\'{n}ski~\cite{DL22} recently provided three possible candidates partially solving the conjecture.
The conjecture in~\cite{GJV05} was based upon the polynomiality of the double Hurwitz numbers.
Moreover, the coefficients of the terms in the (weighted) double Hurwitz number $\mathsf{H}_{d,d}^{g}(\beta)/t!d$
were expressed as the coefficients of certain monomials~\cite[Cor.~$3.8$]{GJV05},
while a few instances were really explicitly determined.

Along this line, it is natural to ask if our quasi-triple Hurwitz
numbers share the polynomiality.
For general quasi-triple Hurwitz numbers, we do not know the answer for the moment.
However,
thanks to the explicit formula of Theorem~\ref{thm:recur2},
the one-part quasi-triple Hurwitz numbers are indeed polynomials.
To present our result, we define speculated intersection numbers
using Witten-like symbols as follows:
\begin{align}
	\langle\langle \tau_{b_1} \tau_{b_2} \cdots \tau_{b_n} \rangle\rangle_g^{d,m}= [\beta_1^{b_1} \beta_2^{b_2}\cdots \beta_n^{b_n}] 
	\frac{(-1)^{b_1+b_2+\cdots+ b_n +n} \mathsf{H}_{d,m}^{g}(\beta)}{\frac{(d-1)!  d^t t!}{(t+1+d)!}   }.
\end{align}

\begin{theorem}[Polynomiality]\label{thm:coeffi}
	For fixed $g,d$ and $m$, the one-part quasi-triple Hurwitz number $\mathsf{H}_{d,m}^{g}(\beta)$ is a symmetric polynomial in $\beta_1, \beta_2, \ldots, \beta_n$ with the highest degree at most $2g$.
	
	Moreover, let $t=2g+n-1+m-d$ and  $b=b_1+b_2+\cdots + b_n$. Then,
	\begin{align}
	\langle\langle \tau_{b_1} \tau_{b_2} \cdots \tau_{b_n} \rangle\rangle_g^{d,m} =  {b_1+b_2+\cdots+ b_n +n \choose b_1+1, b_2+1, \ldots, b_n+1} \sum_{k\geq 0}  (-1)^k \stirlingI{m+k}{m} R_{m+k}(b_1,\ldots, b_n),
\end{align}
	where
 \begin{align}\label{eq:coeff-gen}
 	{R_i(b_1,\ldots,b_n)}
 	&= 	\sum_{j=0}^{d-i} \, \sum_{ \{z_1,\ldots,z_p\} \atop p\geq 0} \, \sum_{k\geq b+n-p -b_{z_1}-\cdots -b_{z_{p}}=q}  (-1)^{p+j+k} \left(\frac{d-1}{2}\right)^{k-q} B_{t+j+1-k}^{(j+1)} \nonumber
 	\\
 	&  \times   {b+n\choose b+n-q}^{-1} \stirlingI{d-i-j}{b+n-q} {t+1+d \choose i, \, d-i-j, \, t+j+1-k, \, q, \, k-q}   .
 \end{align}
\end{theorem}

As a corollary, we have the following surprising result.
\begin{corollary}[Genus zero formula]\label{cor:gen-0}
	The quasi-triple Hurwitz number $\mathsf{H}_{d,m}^{0}(\beta)$ does not depend on the exact values of $\beta_i$ for $1\leq i \leq n$,
	and is given by
	\begin{multline}
\frac{ (d-1)!  d^t t!}{(-1)^n m! }  	\sum_{j=0}^{d-m} \sum_{p=0}^n \,\sum_{k\geq n-p}  
 {(-1)^{p+j+k} \stirlingI{d-m-j}{p}} 
 {n \choose t+j+1} \\
   {t+j+1\choose k}
{k \choose  n-p} \left(\frac{d-1}{2}\right)^{k-n+p} B_{t+j+1-k}^{(j+1)} , \qquad \qquad
	\end{multline}
{
	where $t=n-1+m-d$. In particular, the genus zero one-part double Hurwitz number
	\begin{align}
		\mathsf{H}_{d,d}^{0}(\beta)= t! d^{t-1}=(n-1)! d^{n-2}.
	\end{align}
}
\end{corollary}

{

	 Since the established polynomiality of single Hurwitz numbers by the ELSV formula, the polynomiality
	for other classes of Hurwitz numbers (with at most two non-simple branch points) has been widely discussed. Specially, it is well understood that,
	in general, double Hurwitz numbers admit a piecewise polynomiality even if one fixes the ramification profile for one of the two
	nonsimple branch points. But, if the fixed profile corresponds to a partition with one part, the corresponding Hurwitz numbers is really given by
	a polynomial in the parts of the other nonsimple branch points.
	
	The polynomiality in Theorem 3.5 may be viewed as a kind of graded or piecewise polynomiality since the cover degree $d$
	has to be fixed in advance, and
	our formulas for the studied quasi-triple Hurwitz numbers suggest that these numbers is less likely a polynomial.
However	in the following, we will show that certain class of quasi-triple Hurwitz numbers is surprisingly a polynomial,
generalizing the polynomiality of one-part double Hurwitz numbers.

	\begin{theorem}\label{poly-general}
		For any fixed integer $q \geq 0$, $n>0$ and $g\geq 0$ such that $2g-q+n \geq 3$, the scaled Hurwitz number $\mathsf{H}_{|\beta|,|\beta|-q}^{g}(\beta) /t! |\beta|$ with $\beta=(\beta_1, \beta_2, \ldots, \beta_n)$ is a polynomial in $\beta_1, \beta_2, \ldots, \beta_n$ of the lowest degree at least $2g+n-2-\delta_{0,q}-q$ and the highest degree exactly $4g+n-3$.
	\end{theorem}
	
	\proof  We first have from the proof of Theorem 3.3 that
	\begin{align*}
		& \mathbb{M}:= [y^{2g-q}]
		\bigg\{ x^{\frac{d-1}{2}}   \cdot \prod_{i=1}^d \Big(\frac{\sinh(\frac{i}{2}\log x+\frac{i}{2}y)}{iy/2} \Big)^{c_i}
		\bigg\}   \\
		=& [y^{2g-q+n-1}] \frac{1}{\prod_i \beta_i}  e^{y(d-1)/2} \frac{1}{ x- e^{-y}} \prod_{i=1}^n (x^{\beta_i} - e^{- \beta_i y}) \\
		=& [y^{2g-q +n -1}] \frac{1}{\prod_i \beta_i}  e^{y(d-1)/2} [x^{\beta_1-1} + x^{\beta_1-2} e^{-y} + \cdots + e^{-y(\beta_1-1)}] \prod_{i=2}^n (x^{\beta_i} - e^{- \beta_i y}).
	\end{align*}
	In terms of the power series of $y$, the last formula equals $\frac{1}{\prod_i \beta_i} $ times
	\begin{align*}
		&[y^{2g-q+n-1}]  \bigg\{ \sum_{k\geq 0} (\frac{d-1}{2})^k \frac{y^k}{k!} \bigg\} \bigg\{ \sum_{k\geq 0} \frac{y^k}{k!} \sum_{j=0}^{\beta_1-1} x^j [-(\beta_1 -1-j)]^k  \bigg\} 
		\bigg\{ \prod_{i=2}^n \Big( x^{\beta_i} - \sum_{k \geq 0} \frac{(-\beta_i)^k}{k!} y^k  \Big) \bigg\}  \\
		& = \sum_{l_1 \geq 0, \, l_2 \geq 0} \frac{(\frac{d-1}{2})^{\l_2} }{l_2!}  \times  \frac{ \sum_{j=0}^{\beta_1-1} x^j [-(\beta_1 -1-j)]^{l_1} }{l_1 !}   \\
		&\qquad \times \sum_{2\leq k_1 <\cdots < k_p \leq n, \, p\geq 0, \atop k_1'+\cdots + k_p' =2g-q+n-1-l_1 -l_2, \, k_i'>0} (-1)^p \frac{(-\beta_{k_1})^{k_1'}}{k_1'!} \cdots \frac{(-\beta_{k_p})^{k_p'}}{k_p'!} \prod_{j \in [n] \setminus \{1,k_1, \ldots, k_p\}} (x^{\beta_j} -1) .
	\end{align*}
	Writing the last formula as a power series of $x-1$, we obtain
	\begin{align*}
		\mathbb{M} &= \frac{1}{\prod_i \beta_i} \sum_{l_1 \geq 0, \, l_2 \geq 0} \frac{(\frac{d-1}{2})^{\l_2} }{l_2!} \sum_{i \geq 0}\frac{ \sum_{j=0}^{\beta_1-1} {j\choose i} [-(\beta_1 -1-j)]^{l_1} }{l_1 !}  (x-1)^i  \\
		&\quad \times \sum_{2\leq k_1 <\cdots < k_p \leq n, \, p\geq 0, \atop k_1'+\cdots + k_p' =2g-q+n-1-l_1 -l_2, \, k_i'>0} (-1)^p \frac{(-\beta_{k_1})^{k_1'}}{k_1'!} \cdots \frac{(-\beta_{k_p})^{k_p'}}{k_p'!} \prod_{j\neq 1, \, j\neq k_i, \atop 1\leq i \leq p} \Big\{ \sum_{l>0} {\beta_{j} \choose l} (x-1)^l \Big\}
	\end{align*}
	
	{\em Claim~$1$.} For any fixed $i\geq 0$ and $l_1\geq 0$, we have
	\begin{align}
		\sum_{j=0}^{\beta_1-1} {j\choose i} [-(\beta_1 -1-j)]^{l_1} &= (-1)^{l_1} \sum_{j \geq 0} j! \stirlingII{l_1}{j} {\beta_1 \choose i+j+1} \\
		&= \sum_{k\geq 1} \sum_{j \geq 0}  \frac{(-1)^{l_1+i+j+1-k} j! }{(i+j+1)!} \stirlingII{l_1}{j}  \stirlingI{i+j+1}{k} \beta_1^k 
	\end{align}
	which is a polynomial in $\beta_1$ of the lowest degree $1$ and the highest degree $l_1+i+1$.

	First, we have
	\begin{align*}
		\sum_{k \geq 0} k^l z^k &=\sum_{k\geq 0} \sum_{j = 0}^k {k \choose j} j! \stirlingII{l}{j} z^k\\
		&= \sum_{j\geq 0} j! \stirlingII{l}{j}\sum_{k\geq 0} {k\choose j} z^k = \sum_{j\geq 0} j! \stirlingII{l}{j} \frac{z^j}{(1-z)^{j+1}} .
	\end{align*}
	Note that $\sum_{k\geq 0} {k \choose i} z^k = \frac{z^i}{(1-z)^{i+1}}$.
	Consequently, we obtain
	\begin{align*}
		\sum_{j=0}^{\beta_1-1} {j\choose i} [-(\beta_1 -1-j)]^{l_1} &=[z^{\beta_1}] (-1)^{l_1} \sum_{j\geq 0} j! \stirlingII{l_1}{j} \frac{z^{j+i+1}}{(1-z)^{j+2+i}} \\
		& = [z^{\beta_1}] (-1)^{l_1} \sum_{j\geq 0} j! \stirlingII{l_1}{j} \sum_{r\geq 0} {r+i+j+1 \choose r} z^{r+i+j+1} \\
		&= (-1)^{l_1} \sum_{j \geq 0} j! \stirlingII{l_1}{j} {\beta_1 \choose i+j+1}.
	\end{align*}
	Using the basic identity $(\beta_1)_{i+j+1}=\sum_k (-1)^{i+j+1-k} \stirlingI{i+j+1}{k} \beta_1^k$
	leads to
	Claim~$1$.
	
	{\em Claim~$2$.} Let $p_i(r)$ be the probability of a uniformly random permutation with $i+r$ elements to
	be a fixed-point free permutation with $r$ cycles for $r>0$, and $p_0(0)=1$. Then, for any fixed $q\geq 0$ and $q \geq i\geq 0$, the number 
	$$
	\stirlingI{d-q+i}{d-q} \frac{d!}{(d-q+i)!} = \sum_{ k \geq 1} \sum_{r \geq 0}  (-1)^{q+r-k} p_i(r) \stirlingI{q+r}{k} d^k
	$$
	which is a polynomial in $d$ of the lowest degree $1- \delta_{0,q}$ and the highest degree $q+i$.
	
	The desired number in Claim~$2$ is clearly $1$ if $q=0$. For $q>0$, by definition, we first have
	$$
	\stirlingI{d-q+i}{d-q} = \sum_{\lambda \vdash d-q+i, \atop \ell(\lambda)= d-q} |C_{\lambda}|.
	$$
	It is not difficult to see that $\lambda=(1+\mu_1, 1+ \mu_2, \ldots, 1+ \mu_r, 1, \ldots, 1)$ where
	$r\geq 0$, $(\mu_1, \ldots, \mu_ r) \vdash i$ and there are $d-q-r$ parts equal to $1$.
	Thus,
	\begin{align*}
		\stirlingI{d-q+i}{d-q}  \frac{d!}{(d-q+i)!} &= \sum_{\mu=[1^{m_1},\ldots, i^{m_i}] \vdash i, \atop m_1+\cdots + m_i =r, \, r \geq 0} \frac{(d-q+i)!}{\prod_{j=1}^i (1+j)^{m_j} m_j! (d-q-r)! }  \frac{d!}{(d-q+i)!} \\
		&=\sum_{\mu=[1^{m_1},\ldots, i^{m_i}] \vdash i, \atop m_1+\cdots + m_i =r, \, r \geq 0} \frac{1}{\prod_{j=1}^i (1+j)^{m_j} m_j!  }  d(d-1)\cdots (d-q-r+1) \\
		&= \sum_{ k \geq 1} \sum_{r \geq 0} p_i(r) (-1)^{q+r-k} \stirlingI{q+r}{k} d^k .
	\end{align*}
	The last formula is clearly a polynomial in $d$ of the lowest degree $1$ and the highest degree $q+i$,
	completing the proof of Claim~$2$.
	
	Next, according to Theorem 3.3, we have
	\begin{align*}
		\frac{\mathsf{H}_{|\beta|,|\beta|-q}^{g}(\beta) } {(2g-q+n-1)! |\beta| } & = d^{2g-q+n-3} \left. \frac{\mathrm{d}^{d}}{ \mathrm{d} x^d} \frac{(\log x)^{d-q}}{(d-q)!} \mathbb{M} \right|_{x=1}  \\
		& = d^{2g-q+n-3} d! \cdot [(x-1)^d] \frac{(\log x)^{d-q}}{(d-q)!} \mathbb{M} \\
		& = d^{2g-q+n-3} d! \cdot \sum_i [(x-1)^i] \frac{(\log x)^{d-q}}{(d-q)!}  [(x-1)^{d-i}] \mathbb{M} .
	\end{align*}
	Summarizing the above discussions, we subsequently obtain
	\begin{align*}
		& \frac{\mathsf{H}_{|\beta|,|\beta|-q}^{g}(\beta) } {(2g-q+n-1)! |\beta| } 
		=  \frac{ d^{2g-q+n-3} }{\prod_l \beta_l}  \sum_{i=0}^q \frac{(-1)^i d ! }{ (d-q+i) ! } \stirlingI{d-q+i}{ d-q}    \sum_{l_1 \geq 0, l_2 \geq 0} \frac{(\frac{d-1}{2})^{l_2} }{l_2!} \\
		& \quad \times  \sum_{2\leq k_1 <\cdots < k_p \leq n, \, p\geq 0, \atop k_1'+\cdots + k_p' =2g-q+n-1-l_1 -l_2, \, k_l'>0} (-1)^p \frac{(-\beta_{k_1})^{k_1'}}{k_1'!} \cdots \frac{(-\beta_{k_p})^{k_p'}}{k_p'!} \Bigg\{ \\
		&\qquad \times \sum_{r = n-1-p}^{q-i}  \sum_{ r_{j_1}+\cdots + r_{j_{n-1-p}} =r, \, r_{j_l} >0 \atop \{j_1, \ldots , j_{n-1-p} \} =[n]\setminus \{1, k_1, \ldots, k_p\}} \prod_l {\beta_{j_l} \choose r_{j_l}} 
		\frac{\sum_{j=0}^{\beta_1-1} {j\choose q-i-r} [-(\beta_1 -1-j)]^{l_1} }{l_1 !} \Bigg\} .
	\end{align*}
	For fixed integers $g, q,n$, the appearing indices $l_1, l_2, k_l',r, r_l$ in the last formula are all finite and independent of $\beta_i$.
	Moreover, it is observed that the power of $\beta_i$ for any $i \in [n]$ in the last two lines is at least one.
	Thus, we conclude that the last formula is indeed a polynomial in $\beta_i$'s.
	
	It remains to determine the lowest and highest degrees of the polynomial.
	For the possible highest degree, we can sum up the contributed highest degrees of the factors:
	\begin{align*}
		& (2g-q-3)+(q+i) + (l_2 )+ (2g-q+n-1 -l_1 -l_2) +(r)+(l_1+ q-i-r +1) \\
		= & 4g+n-3 .
	\end{align*}
	Next, we consider the coefficient of the term $\beta_1^{4g+n-3}$ which is given by
	\begin{align*}
		& \sum_{i=0}^q \frac{(-1)^i}{2^i i!}   \sum_{ l_2 \geq 0} \frac{1 }{ 2^{l_2} l_2!} 
		\times   \sum_{l_1\geq 0, \, 0 \leq 2g-i-l_1-l_2 \leq q-i }  \sum_{2\leq k_1 <\cdots < k_p \leq n, \atop  p =2g-q+n-1-l_1 -l_2}
		\frac{(-1)^{l_1} }{(2g-i-l_2+1)! }   \\ 
		= &  \sum_{i=0}^q \frac{(-1)^i}{2^i i!}   \sum_{ l_2 \geq 0} \frac{1 }{ 2^{l_2} l_2!} 
		\times   \sum_{l_1\geq 0, \, 0 \leq 2g-i-l_1-l_2 \leq q-i } {n-1 \choose q+l_1 +l_2-2g}  
		\frac{(-1)^{l_1} }{(2g-i-l_2+1)! } \\
		=& \sum_{i=0}^q \frac{(-1)^i}{2^i i!}   \sum_{ l_2 \geq 0} \frac{1 }{ 2^{l_2} l_2!} 
		\times  {n-2 \choose q-i}  
		\frac{(-1)^{l_2 -i } }{(2g-i-l_2+1)! } \\
		=& \sum_{i=0}^q \frac{1}{2^i i!}  {n-2 \choose q-i}     
		\frac{1 }{2^{2g-i+1} (2g-i+1)! } >0 .
	\end{align*}
	In the above derivation, the following well-known partial sum is used:
	\begin{align}
		\sum_{k=0}^M (-1)^k {N \choose k}=(-1)^M {N-1 \choose M}.
	\end{align}
	As a result, we conclude that the highest degree is exactly $4g+n-3$.

	Analogously, we obtain the lowest possible degree as
	\begin{align*}
		(2g-q-3)+(1-\delta_{0,q}) +( 0) + (p)  +(n-1-p) +1 
		= & 2g+n-3 +1- \delta_{0, q} -q .
	\end{align*}
	This completes the proof. \qed

}

Setting $q=0$ in Theorem 3.6 immediately results in the polynomiality of one-part double Hurwitz numbers.
We can derive all the coefficients of one-part double Hurwitz numbers from Theorem 3.6.

When $m=d$ in Theorem 3.4, we obtain the following stronger result regarding double Hurwitz numbers, where we do not need to fix $d$ in advance
and it is implicitly determined by $\beta$.
First, in order for aligning with the discussion in Goulden, Jackson and Vakil~\cite{GJV05},
we adopt a slightly different version of the Witten symbols.
For $b_1+b_2 + \cdots + b_n +2k=4g-3+n$ where $b_i \geq 0$, $k\geq 0$, and $(g,n) \neq (0,1), (0,2)$, let 
$$
\langle\langle  \tau_{b_1} \tau_{b_2} \cdots \tau_{b_n} \Lambda_{2k}\rangle\rangle_{g}' := (-1)^k [\beta_1^{b_1} \beta_2^{b_2} \cdots \beta_n^{b_n}] \frac{\mathsf{H}_{d,d}^{g}((\beta_1,\ldots, \beta_n))}{t! d}.
$$
For other situations, $\langle\langle  \tau_{b_1} \tau_{b_2} \cdots \tau_{b_n} \Lambda_{2k}\rangle\rangle_{g}' :=0$.  In~\cite{GJV05}, these conjectural intersection numbers were studied.
In view of Corollary~\ref{cor:gen-0}, the case $g=0$ is trivial. In the following, we
assume $g>0$ unless explicitly stated otherwise.

\begin{lemma} The following holds:
\begin{align}
{\bf C}_{b,g,n} &:= \frac{1}{(b+3-2g)!}\sum_{k=0}^{2g+n}   \frac{1}{(k+2g-b-3)!}  \frac{B_{2g+n-k}}{(2g+n-k)!2^k} \nonumber \\
&=\frac{1-2^{4g+n-b-4}}{ 2^{2g+n-1} (b+3-2g)!} \frac{ B_{4g+n-b-3}}{(4g+n-b-3)!}.
\end{align}
\end{lemma}
\proof Recall that $\sum_{n \geq 0} \frac{B_n}{n!} x^n = \frac{x}{e^x-1}$.
Then, we have
\begin{align*}
\frac{x}{e^x-1} e^{x/2} = \sum_{n\geq 0} \sum_{k \geq 0} \frac{B_{k}}{k!} \frac{(1/2)^{n-k}}{(n-k)!} z^n.
\end{align*}
On the other hand, we easily observe
\begin{align}
\frac{x}{e^x-1} e^{x/2} = 2 \frac{x/2}{e^{x/2} -1} -\frac{x}{e^x-1}= \sum_{n\geq 0} \frac{(2^{1-n}-1)B_n}{n!} x^n.
\end{align}
Consequently, we obtain
\begin{align}
\frac{(2^{1-n}-1)B_n}{n!} = \sum_{k \geq 0} \frac{B_{k}}{k!} \frac{(1/2)^{n-k}}{(n-k)!}.
\end{align}
For $b+3-2g \geq 0$, we have 
$$
\sum_{k \geq 0} \frac{B_{k}}{k!} \frac{(1/2)^{4g+n-b-3-k}}{(4g+n-b-3-k)!} = \sum_{k = 0}^{2g+n} \frac{B_{k}}{k!} \frac{(1/2)^{4g+n-b-3-k}}{(4g+n-b-3-k)!}
=\frac{(2^{b+4-4g-n}-1)B_{4g+n-b-3}}{(4g+n-b-3)!},
$$
from which we obtain the closed formula for ${\bf C}_{b,g,n}$ in the lemma. \qed

\begin{theorem}\label{thm:double-coeff}
For	fixed $g$ and $n$,
the number $\mathsf{H}_{d,d}^{g}((\beta_1,\ldots, \beta_n))/ (t! d)$ is a symmetric polynomial of the highest degree $4g-3+n$
and lowest degree $2g+n-3$.
The coefficient of the term $\beta_1^{b_1}\cdots \beta_n^{b_n}$ with $b=b_1+\cdots+b_n$ and $b+2k=4g+n-3$ is given by
\begin{align}
&(-1)^{\frac{4g-3+n-b}{2}}  \langle\langle  \tau_{b_1} \tau_{b_2} \cdots \tau_{b_n} \Lambda_{2k}\rangle\rangle_{g}'
=  {{\bf C}_{b,g,n}} \, \widehat{{\bf C}}_g(b_1,\ldots, b_n)	,
\end{align}
where 
\begin{align}
&\widehat{{\bf C}}_g(b_1,\ldots, b_n)=[x^{2g+n-3}] (2g+n-3)! (b+3-2g)!  \prod_{i=1}^n \frac{(x+1)^{b_i+1} - (x-1)^{b_i+1}}{(b_i+1)!}.
\end{align}
In particular, if $b+n$ is even, then
$$
\langle\langle  \tau_{b_1} \tau_{b_2} \cdots \tau_{b_n} \Lambda_{2k}\rangle\rangle_{g}'=0.
$$
Moreover, if $b_1+\cdots + b_n =2g+n-3$, then
\begin{align}\label{eq:lambda-g}
\langle\langle  \tau_{b_1} \tau_{b_2} \cdots \tau_{b_n} \Lambda_{2g}\rangle\rangle_{g}' = (-1)^g {2g+n-3 \choose b_1, \ldots, b_n} \frac{1-2^{2g-1}}{2^{2g-1} (2g)!} B_{2g}.
\end{align}
\end{theorem}

\proof Based on eq.~\eqref{eq:H-W}, when $m=d$, the coefficient of the term $\beta_1^{b_1}\cdots \beta_n^{b_n}$
equals the one in the following expression:
\begin{align*}
	 &\frac{(\beta_1+\cdots + \beta_n)^{t-2} }{\prod_i \beta_i} 
	\sum_{z_1<\cdots < z_p, \atop p\geq 0} \frac{(-1)^{p}}{(t+1)!} {\bf B}_{t+1} \bigg(\frac{1-\beta_{z_1}-\cdots -\beta_{z_p}+\beta_{z'_1}+\cdots+\beta_{z'_{n-p}}}{2}\bigg).
\end{align*}
The latter is easily seen to be
\begin{align*}
	&\sum_{a_i\leq b_i+1,\atop
	1\leq i \leq n} {t-2 \choose a_1,\ldots, a_n} \sum_{z_1<\cdots < z_p, \atop p\geq 0} \frac{(-1)^{p}}{(t+1)!}
\sum_{k=0}^{t+1} {t+1\choose k} \frac{1}{2^k} {k\choose b+n-t+2}\\
& \qquad \qquad \times (-1)^{p+\sum_i b_{z_i}-a_{z_i}} {n+b-t+2 \choose b_1+1-a_1,\ldots, b_n+1-a_n} B_{t+1-k},\\
&= {\bf C}_{b,g,n}	\sum_{a_i\leq b_i+1,\atop
	1\leq i \leq n} {2g+n-3 \choose a_1,\ldots, a_n}  {b+3-2g \choose b_1+1-a_1,\ldots, b_n+1-a_n}  \sum_{z_1<\cdots < z_p, \atop p\geq 0} 
(-1)^{\sum_i b_{z_i}-a_{z_i}}.
\end{align*}
Note that for a fixed $(a_1,\ldots, a_n)$ with $a_1+\cdots + a_n =2g+n-3$,
if $\sum_{i=1}^n b_i-a_i =b-t+1$ is odd, i.e., $b+n$ is even, then $(-1)^{\sum_{i=1}^p b_{z_i}-a_{z_i}}=-(-1)^{\sum_{q\notin \{z_1,\ldots, z_p\}} b_q-a_q}$.
Consequently, their contributions to the coefficient cancel each other and the coefficient in this case is thus zero.
Next, we observe that the ``big" factor in the last formula also equals
$$
[x^{2g+n-3}] (2g+n-3)! (b+3-2g)!  \prod_{i=1}^n \frac{(x+1)^{b_i+1} - (x-1)^{b_i+1}}{(b_i+1)!}.
$$
Note there is a sign $(-1)^k= (-1)^{\frac{4g-3+n-b}{2}}$ for $\langle\langle  \tau_{b_1} \tau_{b_2} \cdots \tau_{b_n} \Lambda_{2k}\rangle\rangle_{g}'$.

Obviously, ${{\bf C}_{b,g,n}}=0$ if $b>4g+n-3$. Thus, the highest degree is at most $4g+n-3$.
To verify the sharpness, it suffices to consider the coefficient of the term $\beta_1^{4g+n-3} \beta_2^0 \cdots \beta_{n}^0$
which is easily seen to be non-zero from our expression.
In the other end, it is also obvious that $\widehat{{\bf C}}_g(b_1,\ldots, b_n)=0$ if $b < 2g+n-3$.
Now suppose $b=2g+n-3$. Note that the highest degree term contributed by $(x+1)^{b_i+1}- (x-1)^{b_i+1}$ is exactly $2 (b_i+1) x^{b_i}$.
As a result, it is clear that
$$
[x^{2g+n-3}] (2g+n-3)! (b+3-2g)!  \prod_{i=1}^n \frac{(x+1)^{b_i+1} - (x-1)^{b_i+1}}{(b_i+1)!} = (2g+n-3)! n! \cdot 2^n \prod_i \frac{1}{b_i !}.
$$
Considering that
$$
{{\bf C}_{2g+n-3,g,n}}= \frac{1}{2^n n!} \frac{ 1- 2^{2g-1}}{2^{2g-1} (2g)!} B_{2g},
$$
we easily arrive at eq.~\eqref{eq:lambda-g}.
This completes the proof.
\qed

The special case, eq.~\eqref{eq:lambda-g}, was obtained in~\cite{GJV05} and is an analogue of the celebrated $\lambda_g$-conjecture (proved in~\cite{FP}).
As an application of Theorem~\ref{thm:double-coeff}, we may provide a new proof of the string and dilaton equations proved in~\cite{GJV05}.
The proof is left to the interested reader.

\begin{proposition}[String and dilaton equations~\cite{GJV05}]
  The following equations are true:
  \begin{align}
\langle\langle  \tau_0 \tau_{b_1} \tau_{b_2} \cdots \tau_{b_n} \Lambda_{2k}\rangle\rangle_{g}'
&= \sum_{i=1}^n 
\langle\langle  \tau_{b_1}\cdots \tau_{b_{i-1}} \tau_{b_i-1}  \tau_{b_{i+1}}\cdots \tau_{b_n} \Lambda_{2k}\rangle\rangle_{g}' ,\\
\langle\langle \tau_1 \tau_{b_1} \tau_{b_2} \cdots \tau_{b_n} \Lambda_{2k}\rangle\rangle_{g}' &=(2g+n-2) \langle\langle  \tau_{b_1} \tau_{b_2} \cdots \tau_{b_n} \Lambda_{2k}\rangle\rangle_{g}'.
\end{align}
\end{proposition}

 As mentioned, only a few coefficients of the one-part double Hurwitz numbers were really explicitly determined before.
On the contrary, instead of knowing the coefficients of certain terms, we have the advantage
of knowing the coefficients of any terms according to Theorem~\ref{thm:double-coeff} and Theorem~\ref{thm:coeffi}.
This may facilitate the search of ELSV-type expressions for the one-part double
and quasi-triple Hurwitz numbers.
These related questions will be studied in our future work.

\begin{remark}
	Those who are experts on Bernoulli polynomials and numbers may be able to further simply some of our formulas in this section,
	and we are looking forward to that.
\end{remark}

\section{Future studies}

This is the first of a series of work on studying Hurwitz numbers and hypermaps as well as constellations.
We presented the study of Hurwitz numbers mainly from a combinatorial point of view.
Due to the rich connections of single and double Hurwitz numbers to enumerative geometry
and mathematical physics, in general, there is no doubt that triple Hurwitz numbers are also closely related.
But, triple Hurwitz numbers are even harder to compute and not much information has been obtained for them before.
Thus, almost no concrete connections have been established.
Hopefully, the presented structure in this work may help in finding concrete connections.
This will be our first direction of future exploration.

Secondly, we are just at the beginning of exploring the applications of the novel recursion.
We will apply it to other problems admitting a permutation product model.
Some works along this line are in progress, 
for instance, the study of hypermaps and maps~\cite{chen-1,b-chen}, as well as Hurwitz numbers with completed cycles.

As for Hurwitz numbers with completed cycles, Okounkov and Pandharipande~\cite{ok-pa} showed they
are equivalent to the stationary sector of the relative Gromov-Witten theory.
Thus, one may compute the latter via computing the former. In~\cite{ok-pa}, certain
generating function for double Hurwitz numbers with completed cycles was expressed in the operator formalism
on the infinite wedge space. Double Hurwitz numbers with completed cycles all being a fixed completed $r$-cycle
were also studied using the operator formalism on the infinite wedge space~\cite{complete11}, where in particular,
the piece-wise polynomiality of the double Hurwitz numbers with only completed $r$-cycles (for a fixed $r$) were obtained.
In Nguyen~\cite{3-cycle}, an explicit formula for one-part double Hurwitz numbers with only completed $3$-cycles
was obtained, and the polynomiality for this special case follows consequently.
We can somehow generalize these works by giving explicit formulas for double and (quasi) triple Hurwitz
numbers with arbitrary combination of completed cycles (instead of the same fixed $r$-cycle), and also deriving
their polynomiality. These results and related will be reported in the work~\cite{chw} soon.


\section*{Acknowledgements}

The author thanks Prof.~Sergey Shadrin for pointing out that Siegel-Veech weight may be related to triple Hurwitz numbers.


\end{document}